\documentclass[12pt,twoside,a4paper]{article}
\usepackage{amsmath,amssymb,latexsym,theorem,bbm,natbib,epsfig,color}
\usepackage{multirow,graphicx,array,epstopdf,easybmat}  
\newfont{\msa}{msam10 scaled\magstep1}
\newfont{\ssmsa}{msam9}

\newcommand{\tr}{\mathrm{tr}}
\newcommand{\mspe}{\mathrm{MSPE}}
\newcommand{\imspe}{\mathrm{IMSPE}}
\newcommand{\ent}{\mathrm{Ent}}

\numberwithin{equation}{section}

\newcommand{\proofend}{\hfill$\square$}

\newtheorem{thm}{Theorem}[section]

\theorembodyfont{\rm}

\newtheorem{ex}[thm]{Example}

\title{On the optimal designs for the prediction of complex Ornstein-Uhlenbeck processes}

\author{ {\sc Kinga Sikolya} and {\sc S\'andor Baran}\\ 
         Faculty of Informatics, University of Debrecen\\
         Kassai \'ut 26, H-4028 Debrecen, Hungary}

\date{}

\begin{document}
\pagestyle{myheadings}

\maketitle

\begin{abstract}
Physics, chemistry, biology or finance are just some examples out of the many fields where complex Ornstein-Uhlenbeck (OU) processes have various applications in statistical modelling. They play role e.g. in the description of the motion of a charged test particle in a constant magnetic field or in the study of rotating waves in time-dependent reaction diffusion systems, whereas Kolmogorov used such a process to model the so-called Chandler wobble, the small deviation in the Earth's axis of rotation. A common problem in these applications is
deciding how to choose a set of a sample locations in order to predict a random process in an optimal way.
We study the optimal design problem for the prediction of a complex OU process on a compact interval with respect to integrated mean square prediction error (IMSPE) and entropy criteria. We derive the exact forms of both criteria, moreover, we show that optimal designs based on entropy criterion are equidistant, whereas the IMSPE based ones may differ from it. Finally, we present some numerical experiments to illustrate selected cases of optimal designs for small number of sampling locations.

\end{abstract}

\section{Introduction}
  \label{sec:sec1}
Physics, chemistry, biology or finance are just some examples out of the many fields where random processes have various applications in statistical modelling. In the current paper we study the problem of optimal design for the prediction of a complex Ornstein-Uhlenbeck (OU) process on a compact interval with respect to the Integrated Mean Square Prediction Error (IMSPE) and entropy criteria. A complex OU process \citep[see e.g.][]{arato}, defined in detail in Section \ref{sec:sec2}, is used in the study of rotating waves in time-dependent reaction diffusion systems \citep{bl,otten}, it can also describe the motion of a charged test particle in a constant magnetic field \citep{balescu}, and has several applications in financial mathematics as well \citep[see e.g.][]{bns}. An important application of the complex OU
is the so-called Chandler wobble, the small deviation in the Earth's axis of rotation relative to the solid earth \citep{lambeck}.
The uniqueness of this motion is the change of about 9 meters (30 ft) in the point where the axis intersects the Earths surface and it has a period of 435 days (14 months). Recently monitoring of the polar motion is done by the International Earth Rotation Service (IERS).
Kolmogorov described the Chandler Wobble using the model
\begin{equation}
  \label{chandler}
Z(t)=Z_1(t)+iZ_2(t)=m{\mathrm e}^{i2\pi t}+Y(t), \qquad t>0,
\end{equation}
where \ $Z_1(t)$ \ and \ $Z_2(t)$ \ are the coordinates of the deviation of the instantaneous pole from the North Pole \citep{aks}. 
In this model the first term is a periodical component, whereas the second one \ $Y(t)$ \ is a complex OU process.

For the complex OU process \citet{bszs} derived the exact form of the Fisher information matrix (FIM) and investigated the properties of D-optimal designs for estimation of model parameters. Here we derive
exact form of the IMSPE already studied e.g. in \citet{bss13}, \citet{bss15} or \citet{Crary}, and show that in contrast to the case of classical real
OU process on a compact interval investigated by \citet{baz}, the equidistant design is usually not optimal. 
We also investigate the properties of the optimal design with respect to the entropy criterion \citep{sw,baz}.
The paper is organized as follows. Section \ref{sec:sec2} gives the necessary definitions connected to the complex OU process, Sections \ref{sec:sec3} and \ref{sec:sec4} contain the results about IMSPE and entropy criteria, respectively, whereas in Section \ref{sec:sec5} some simulation results are presented. The paper ends with the concluding Section \ref{sec:sec6} and to maintain the continuity of the explanation, proofs are provided in the Appendix.

\section{Complex Ornstein-Uhlenbeck process with a trend}
  \label{sec:sec2}

Let \ $Y(t)=Y_1(t)+iY_2(t), \ t\geq 0$, \ be a complex valued stationary OU process defined by the stochastic differential equation (SDE)
\begin{equation}
  \label{c_sde}
{\mathrm d}Y(t)=-\gamma Y(t){\mathrm d}t+ \sigma {\mathrm d}{\mathcal W}(t), \qquad Y(0)=\xi,
\end{equation}
where \ $\gamma=\lambda -i\omega$ \ with \ $\lambda>0, \ \omega\in {\mathbb R}$, \ $\sigma>0$,  \ ${\mathcal W}(t)={\mathcal W}_1(t)+i{\mathcal W}_2(t), \ t\geq 0$, \ is a standard complex Brownian motion, with \ ${\mathcal W}_1(t)$ \ and \ ${\mathcal W}_2(t)$ \ being independent standard Brownian motions, and \ $\xi=\xi_1+i\xi_2$, \ where \ $\xi_1$ \ and \ $\xi_2$ \ are centered normal random variables that are chosen according to stationarity \citep{arato}.
In the current paper we consider the shifted complex stochastic process \ $Z(t)=Z_1(t)+iZ_2(t)$, \ defined as
\begin{equation}
  \label{c_mod}
Z(t)=m+Y(t), \qquad t\geq 0,
\end{equation}
where \ $m=m_1+im_2, \ m_1,m_2\in{\mathbb R}$, \ and \ $Y(t), \ t\geq 0$, \ is a complex valued stationary OU process, observed in design points taken from the non-negative half-line \ ${\mathbb R_+}$. 

In order to simplify calculations, instead of the complex process \ $Y(t)$ \ one can use the two-dimensional real valued stationary OU process \ $\big( Y_1(t),Y_2(t)\big)^{\top}$ \ defined by the SDE
\begin{equation}
  \label{2d_sde}
  \begin{bmatrix} {\mathrm d}Y_1(t) \\ {\mathrm d}Y_2(t)\end{bmatrix} =
  A
   \begin{bmatrix} Y_1(t) \\ Y_2(t)\end{bmatrix}{\mathrm d}t
   +\sigma  \begin{bmatrix} {\mathrm d}{\mathcal W}_1(t) \\ 
     {\mathrm d}{\mathcal W}_2(t)\end{bmatrix}, \quad \text{where} \quad 
   A:=\begin{bmatrix} -\lambda & -\omega \\ \omega & -\lambda \end{bmatrix}.
\end{equation}
Note that in physics \eqref{2d_sde} this is called A-Langevin equation, see e.g. Balescu \cite{balescu}.
The real and imaginary parts of a complex OU process form a two-dimensional real OU process satisfying \eqref{2d_sde}, and conversely, if
\ $\big( Y_1(t),Y_2(t)\big)^{\top}$ \ satisfies \eqref{2d_sde} then \  $Y_1(t)+iY_2(t)$ \ will be a complex OU process which solves \eqref{c_sde}.
Naturally, \ ${\mathsf E}  Y_1(t)={\mathsf E}  Y_2(t)=0$, \ and the covariance matrix function of the process \ $\big( Y_1(t),Y_2(t)\big)^{\top}$ \ is given by
\begin{equation}
 \label{ou_cov}
{\mathcal R}(\tau):= {\mathsf E}  \begin{bmatrix} Y_1(t+\tau ) \\ Y_2(t+\tau )\end{bmatrix}  \begin{bmatrix} Y_1(t) \\ Y_2(t)\end{bmatrix}^{\top} \!\!= \frac {\sigma^2}{2\lambda} {\mathrm e}^{A\tau}=\frac {\sigma^2}{2\lambda} {\mathrm e}^{-\lambda \tau}\begin{bmatrix} \cos(\omega\tau)  & \sin(\omega\tau) \\ -\sin (\omega\tau) & \cos (\omega\tau) \end{bmatrix}, \quad \tau \geq 0.
\end{equation}
This results in a complex covariance function of the complex OU process \ $Y(t)$ \ defined by \eqref{c_sde} of the form
\begin{equation*}
{\mathcal C}(\tau):={\mathsf E}Y(t+\tau)\overline{Y(t)}=\frac {\sigma^2}{\lambda} {\mathrm e}^{-\lambda \tau}\big( \cos(\omega\tau)-i\sin(\omega\tau)\big ), \qquad \tau\geq 0,
\end{equation*}
behaving like a damped oscillation with frequency \ $\omega$.

In the next steps we will consider the damping parameter \ $\lambda$, \ frequency \ $\omega$ \ and standard deviation \ $\sigma$ \ as known.
Nevertheless a promising direction moving forward could be the examination of models where the above mentioned parameters should also be estimated. Note that the estimation of \ $\sigma$ \ can be done on the basis of a single realization of the complex process, see e.g. \citet[][Chapter 4]{arato}. Now, without loss of generality, one can set the variances of \ $Y_1(t)$ \ and \ $Y_2(t)$ \ to be equal to one, that is \ $\sigma^2/(2\lambda)=1$, \  which reduces \ ${\mathcal R}(\tau)$ \ to a correlation matrix function. In \citet{abi} we can find more results on the maximum-likelihood estimation of the covariance parameters.

\section{Optimal design with respect to the IMSPE criterion}
  \label{sec:sec3}

Assume we observe our complex process \ $Z(t)$ \ at design points \ $0\leq t_1<t_2 < \ldots <t_n$ \ resulting in the $2n$-dimensional real vector \ $Z=\big( Z_1(t_1),Z_2(t_1),Z_1(t_2),Z_2(t_2), \ldots ,Z_1(t_n),Z_2(t_n)\big)^{\top}$\!\!, where
\begin{equation*}
Z_1(t)=m_1+Y_1(t), \qquad Z_2(t)=m_2+Y_2(t).
\end{equation*} 

The main aim of the kriging technique consists of predicting the output
of the investigated process or field on the experimental region, and for any untried location $x \in {\mathcal {X}}$, which in our case lies in a closed interval ${\mathcal X}\subset {\mathbb R}$, the estimation procedure is
focused on the best linear unbiased estimator (BLUE) \ $\widehat Z(x)$ \ of \ $Z(x)$. \ Considering again the two-dimensional vector form of the complex process, the real and imaginary parts \ $\widehat Z_1(x)$ \ and \ $\widehat Z_2(x)$, \ respectively, of the BLUE are given as
\begin{equation*}
  \begin{bmatrix} \widehat Z_1(x) \\ \widehat Z_2(x)\end{bmatrix} =
   \begin{bmatrix} \widehat m_1 \\ \widehat m_2 \end{bmatrix}
	+{\mathcal Q}(x)C^{-1}(n) \Bigg(\begin{bmatrix}   Z_1 \\  Z_2 \end{bmatrix}
	-H(n)^{\top}\begin{bmatrix} \widehat m_1 \\ \widehat m_2 \end{bmatrix} \Bigg),    
\end{equation*}
where \ $(Z_1 \quad  Z_2)^{\top}:=\big( Z_1(t_1),Z_2(t_1),Z_1(t_2),Z_2(t_2), \ldots ,Z_1(t_n),Z_2(t_n)\big)$ \ is the real observation vector, \  $H(n)$ \ is the \ $2 \times 2n$ \ matrix
  \begin{equation*}
    H(n):=\begin{bmatrix} 1 & 0 & 1 & 0 &\cdots &  1 & 0 \\ 0 & 1 & 0 & 1 &\cdots &  0 & 1 \end{bmatrix},
  \end{equation*}
$C(n)$ \ is the \ $2n \times 2n$ \ covariance matrix of the observations, ${\mathcal Q}(x)$ \ is the \ $2  \times 2n$ \ matrix of correlations between \ $Z(x)=\big(Z_1(x),Z_2(x)\big)$ \ and \ $\big\{ \big(Z_1(t_j),Z_2(t_j)\big), \ j=1,2,\ldots ,n\big\}$ \ defined by \
${\mathcal Q}(x)=\big(Q(x,t_1),\ldots
,Q(x,t_n)\big)$ \ with 
\begin{equation*}
 Q(x,t_i):=\frac {\sigma^2}{2\lambda} {\mathrm e}^{-\lambda |x-t_{i}|}\begin{bmatrix} \cos(\omega(x-t_{i}))  & \sin(\omega(x-t_{i})) \\ -\sin (\omega(x-t_{i})) & \cos (\omega(x-t_{i})) \end{bmatrix},
\end{equation*}
and \ $({\widehat m_1},{\widehat m_2})^{\top}$ \ is the generalized least squares estimator of \ $(m_1,m_2)^{\top}$, that is
\begin{equation*}
  \begin{bmatrix} \widehat m_1 \\ \widehat m_2\end{bmatrix} =
	\Big(H(n)C(n)^{-1}H(n)^{\top}\Big)^{-1}H(n)C(n)^{-1}\begin{bmatrix}   Z_1 \\  Z_2 \end{bmatrix}.
\end{equation*}
Thus, a natural criterion is to minimize suitable functionals of the mean squared prediction error (MSPE) given by
\begin{equation}
   \label{mspe}
\mspe\big(\widehat Z(x)\big):= \frac {\sigma^2}{2\lambda} {\mathrm{tr}}\Bigg[{\mathbb I}_2-\big({\mathbb I}_2,\,{\mathcal
  Q}(x)\big)
\left[
      \begin{BMAT}{c.c}{c.c}
      {\mathbb O}_2 & H_n\\
       H_n^{\top} & C(n)
      \end{BMAT}
\right]^{-1}
\big({\mathbb I}_2,\,{\mathcal
  Q}(x)\big)^{\top}\Bigg],
\end{equation}
where \ $\mathbb I_k$ \ and \ ${\mathbb O}_k, \ k\in{\mathbb N}$, \ denote the $k$-dimensional unit and null matrices, respectively. Since the prediction accuracy is often
related to the entire prediction region \ $\mathcal{X}$, \ the
design criterion IMSPE is given by
\begin{equation*}
\imspe\big (\widehat
Z\big):=\frac {2\lambda}{\sigma^2}\int\limits_{\mathcal{X}}\mspe\big (\widehat
Z(x)\big){\mathrm d}x.
\end{equation*}
Without loss of generality we may assume \
$\mathcal{X}=[0,1]$, \ therefore, as extrapolative prediction
is not advisable in kriging, we can set \ $t_1=0$ \ and \ $t_n=1$.

\begin{thm}
 \label{IMSPE}
\ In our setup,
\begin{equation}
  \label{imspe}
 \imspe\big (\widehat Z\big)=2(1-A_n+G(n)^{-1}B_n),
\end{equation}
where
\begin{align*}
G(n)=&\,\!1\!+\!\sum_{\ell=1}^{n-1} g(d_{\ell}), \quad \text{with} \quad g(x)\!:=\!\frac{1\!-\!2{\mathrm e}^{-\lambda x}\cos(\omega x)\!+\!{\mathrm e}^{-2\lambda x}}{1\!-\!{\mathrm e}^{-2\lambda x}}, \ x\!>\!0, \ \ \text{and} \ \ g(0)\!:=\!0, 
\end{align*}
\begin{align}
A_n=&\,\varrho_{n,n}+\sum_{i=1}^{n-1}\frac  
{\varrho_{i,i}-2{\pi_i}\varrho_{i+1,i}+{{\pi_i}^2}\varrho_{i+1,i+1}}{1-\pi_i^2},  \label{an}\\ \label{bn}
B_n=&\,1-2v_{n,n}+\varrho_{n,n} 
-2\sum_{i=1}^{n-1}\frac{(v_{i,i}-{\pi_i}v_{i,i+1})-{\pi_i}(v_{i+1,i}-{\pi_i}v_{i+1,i+1})}{1-{\pi_i}^2}\\ 
&+2\sum_{i=1}^{n-1}\frac{\big(\varrho_{n,i}-\pi_i\varrho_{n,i+1}\big)\big(\cos(\omega(d_i+\ldots + d_{n-1}))-{\pi_i}\cos(\omega(d_{i+1}+\ldots + d_{n-1}))\big)}{1-{\pi_i}^2}\nonumber \\
&+\sum_{i=1}^{n-1} \frac{\big(\varrho_{i,i}-2{\pi_i}\varrho_{i+1,i}
+{\pi_i}^2\varrho_{i+1,i+1}\big)\big(1-2{\pi_i}\cos(\omega d_i)+{\pi_i}^2\big)}{(1-\pi_i^2)^2} \nonumber\\
&+2\sum_{i=2}^{n-1}\sum_{j=1}^{i-1} \big(\varrho_{i,j} - \pi_i\varrho_{i+1,j} -\pi_j \varrho_{i,j+1}+\pi_i\pi_j\varrho_{i+1,j+1}\big) \nonumber\\
&\phantom{2\sum_{i=2}^{n-1}}\times\bigg(
\frac{\cos(\omega(d_j+\ldots + \!d_{i-1}))\!-\!{\pi_i}\cos(\omega(d_j+\ldots + d_{i}))\!-\!{\pi_j}\cos(\omega(d_{j+1}+\ldots + d_{i-1}))}{(1-{\pi_i}^2)(1-{\pi_j}^2)} \nonumber \\
&\phantom{===========================}+\frac{{\pi_i}{\pi_j}\cos(\omega (d_{j+1}+\!\ldots \!+ d_{i}))}
{(1-{\pi_i}^2)(1-{\pi_j}^2)} \bigg), \nonumber
\end{align}
where for \ $i\land j:=\min\{i,j\}, \ i\lor j:=\max\{i,j\}, \ i,j
\in{\mathbb N}$,
\begin{align*}
\varrho_{i,j}:=&\,\frac{1}{2\lambda} \Big(2{\mathrm e}^{-\lambda (d_j+ \cdots + d_{i-1})} - {\mathrm e}^{-\lambda(2d_1+\cdots +2d_{j-1}+d_j+\cdots +d_{i-1})}- {\mathrm e}^{-\lambda(d_j+\cdots +d_{i-1}+2d_i+\cdots +2d_{n-1})}\Big) \\
&+\big(d_j+ \cdots + d_{i-1}\big){\mathrm e}^{-\lambda(d_j+ \cdots + d_{i-1})}, \quad \quad 1\leq j \leq i \leq n,\\[3mm]
v_{i,j}:=&\,\frac {2\lambda}{\lambda^2 +\omega^2}\cos(\omega (d_{i\land j}+\ldots +d_{i \lor j-1})) \\
&+\frac {{\mathrm e}^{-\lambda (d_1+ \cdots + d_{i-1})}}{\lambda^2 +\omega^2}\big(\omega \sin(\omega (d_1+ \cdots + d_{j-1}))-\lambda \cos (\omega (d_1+ \cdots + d_{j-1}))\big) \\
&+\frac {{\mathrm e}^{-\lambda (d_i+ \cdots + d_{n-1})}}{\lambda^2 +\omega^2}\big(\omega \sin(\omega (d_j+ \cdots + d_{n-1}))-\lambda \cos (\omega (d_j+ \cdots + d_{n-1}))\big), 
\end{align*}	
with the empty sum to be defined as zero, and $\pi_i:=\exp(-\lambda d_i)$ \ with \ $d_i:=t_{i+1}-t_i$, \ $i=1,2, \ldots, n-1$.
\end{thm}

\begin{figure}[t!]
\begin{center}
\leavevmode
\hbox{
\epsfig{file=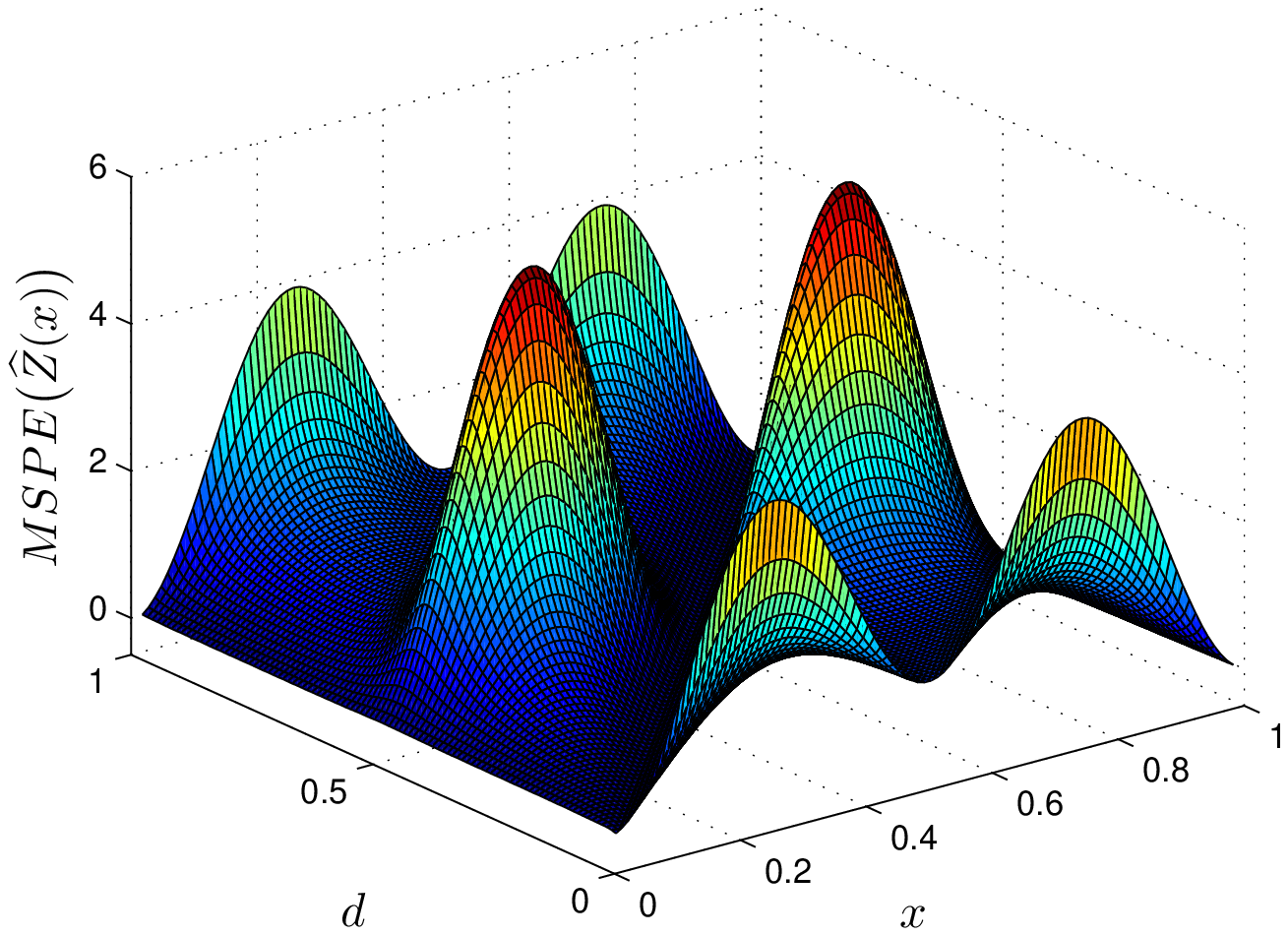,height=6cm, width=7.5 cm} \quad
\epsfig{file=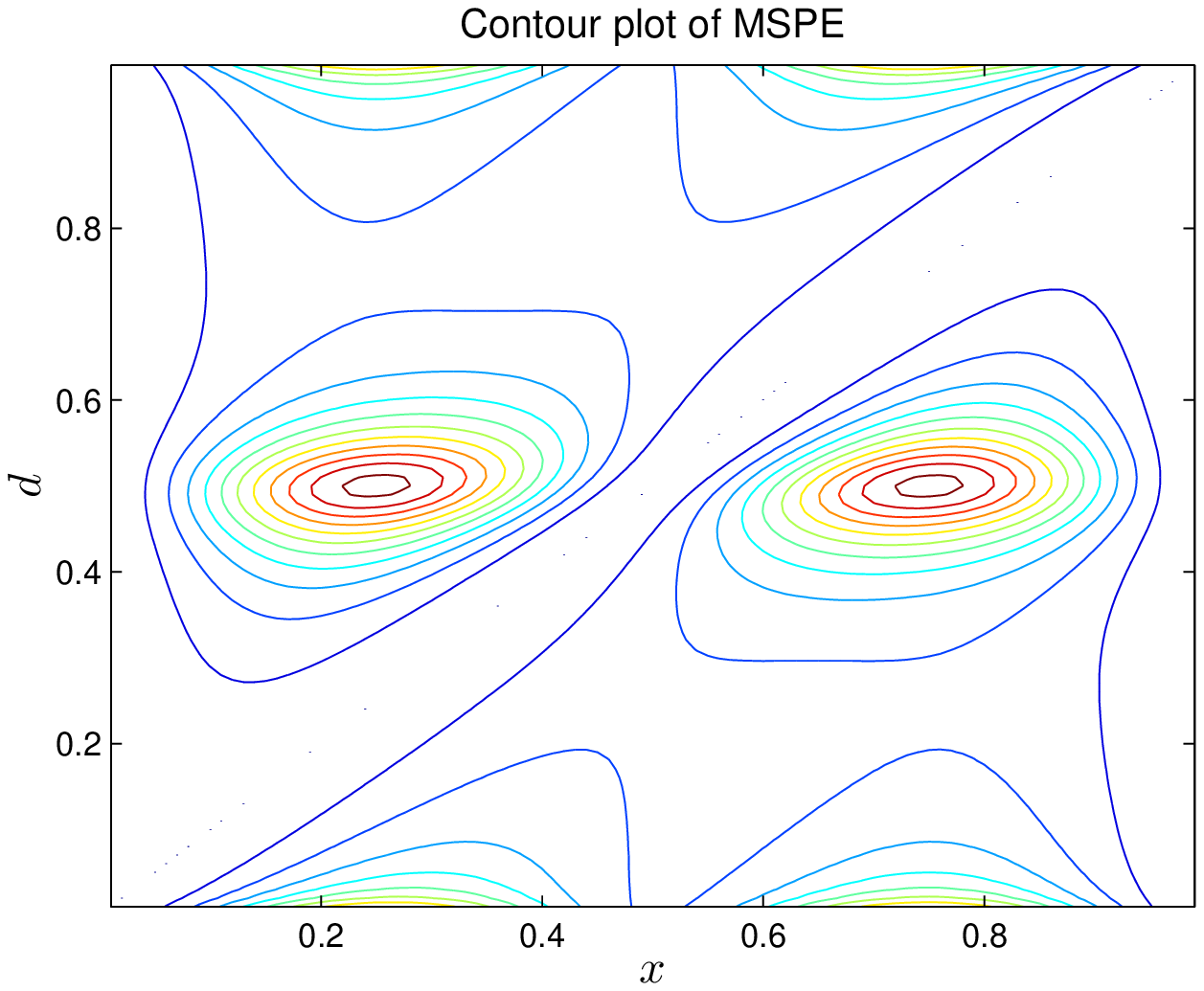,height=5.5cm}}

\centerline{\hbox to 9 truecm {\scriptsize (a) \hfill (b)}}

\smallskip

\leavevmode
\hbox{
\epsfig{file=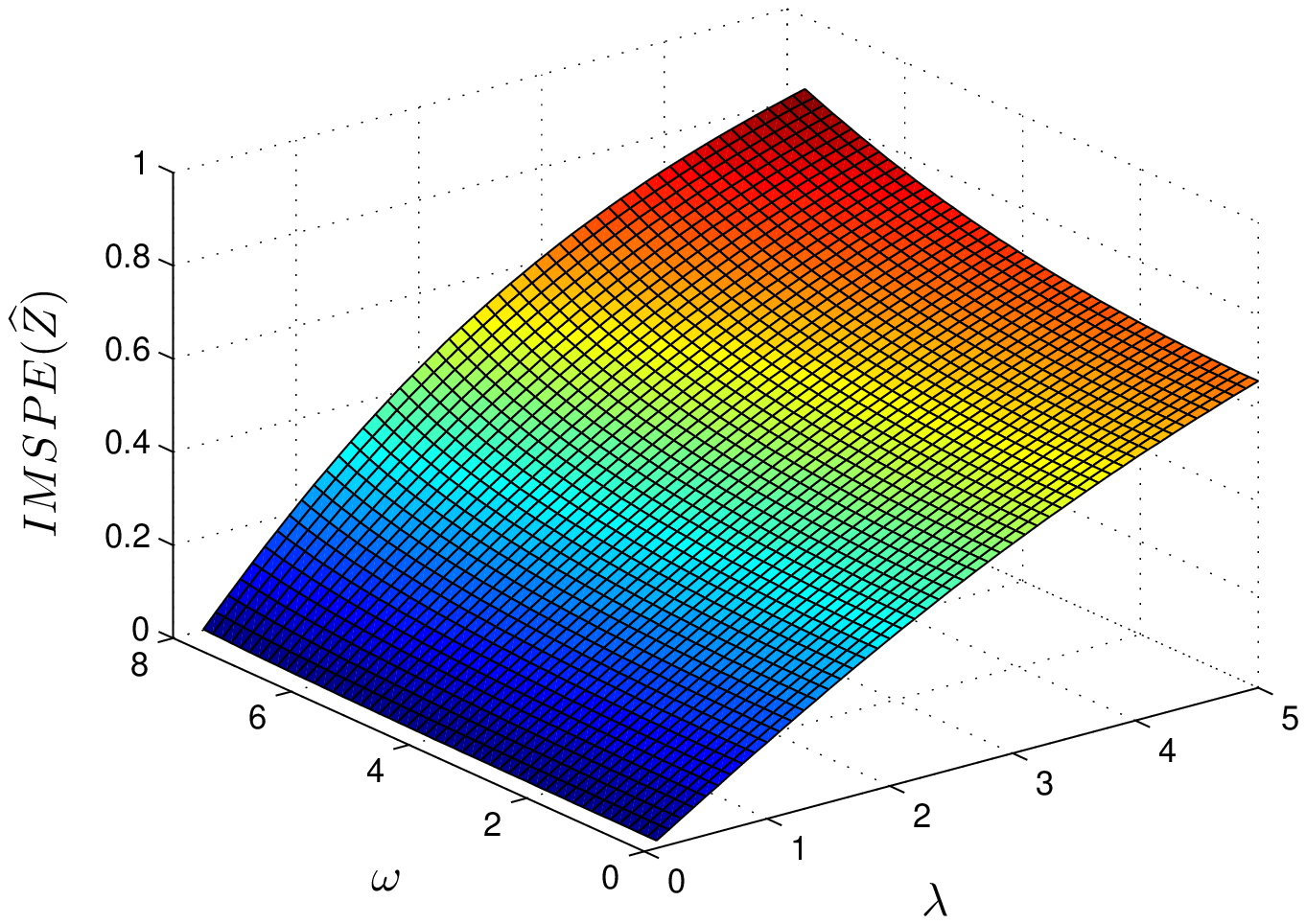,height=6cm, width=7.5 cm} \quad
\epsfig{file=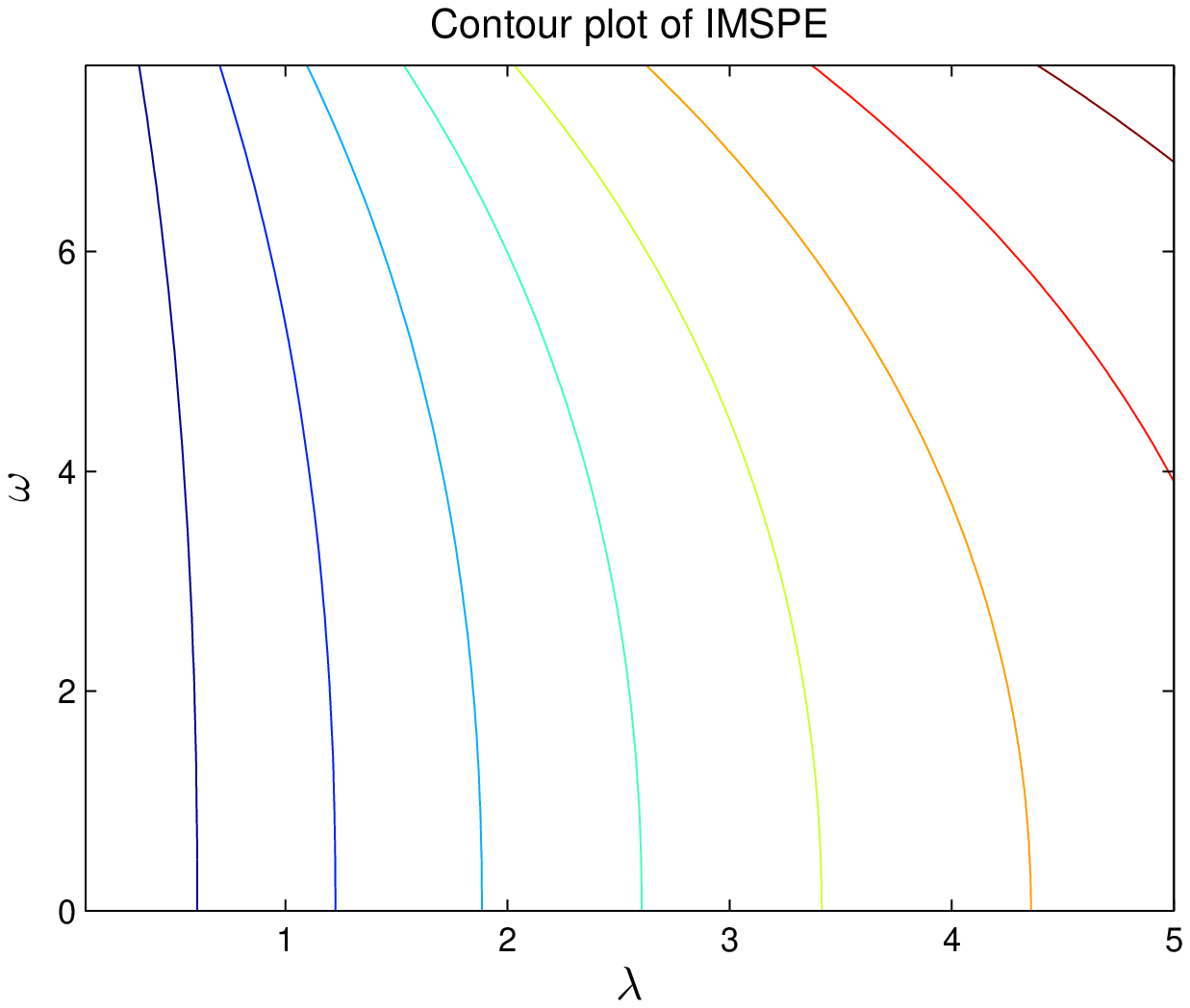,height=5.5cm}}

\centerline{\hbox to 9 truecm {\scriptsize (c) \hfill (d)}}

\end{center}
\caption{Mean squared prediction error (MSPE) for the three-point design \ $\{0,d,1\}$  for \ $\lambda=1, \ \omega=4$ \ and  
IMSPE corresponding to the equidistant three-point design for the prediction as function of \ $\lambda$ \ and \ $\omega$ \ (c), together with the corresponding contour plots (b) and (d), respectively.}
\label{fig1}
\end{figure}

\begin{ex}
  \label{mspe3} \
As an illustration consider a three-point design, that is \ $n=3$, \
$t_1=0, \ t_2:=d, \ t_3=1$.\  Figure \ref{fig1}a shows the mean squared prediction error (MSPE) function
for \ $\lambda=1, \ \omega=4$ \ together with the corresponding contour plot (Figure \ref{fig1}b). 
In Figures \ref{fig1}c and \ref{fig1}d the IMSPE corresponding to the equidistant three-point design for the prediction as function of \ $\lambda$ \ and \ $\omega$ \ and its contour plot, respectively, are given. Straightforward calculation shows that in this case \ $d=\frac{1}{2}$ \ is a minimizer of \ $\imspe\big (\widehat Z\big)$ \ for all possible values of \ $\lambda$ \ and \ $\omega$, \ that is the optimal design is equidistant.
\end{ex}

\begin{ex}
  \label{mspe4} \
Consider now the four-point design \ $\{0,d_1,d_1+d_2,1\}$. \ In this case the partial derivatives of  \ $\imspe\big (\widehat Z\big)$ \ with respect to \ $d_1$ \ and \ $d_2$ \ at \ $d_1=d_2=1/3$ \ are not necessary zero, that is in general, one cannot state that the equidistant design is optimal.
\end{ex}

\section{Optimal information gain for complex OU process}
  \label{sec:sec4}
Another approach to optimal design is to find locations which maximize
the amount of obtained information. Following the ideas of \citet{sw} 
one has to maximize the entropy \ $\ent (Z)$ \ of the
observations corresponding to the chosen design which in our Gaussian
case forms  a $2n$-dimensional normal vector with covariance 
matrix \ $\frac{\sigma ^2}{2\lambda}\,C(n)$, \ that is
\begin{equation*}
\ent (Z)= n\bigg(1+\ln \Big(\frac{\pi \sigma^2}{\lambda}\Big)\bigg)+\frac
12\ln\det C(n).
\end{equation*}

\begin{thm}
  \label{entropy}
Under conditions of Theorem \ref{IMSPE} entropy \ $\ent (Z)$ \
has the form
\begin{equation}
  \label{entropy}
\ent (Z)= n\bigg(1+\ln \Big(\frac{\pi \sigma^2}{\lambda}\Big)\bigg)+\frac
12 \sum_{i=1}^{n-1} \ln \big(1-2\pi_{i}^2 \big).
\end{equation}
For any sample size the equidistant design \ $d_1=d_2=\ldots=d_{n-1}$ \
is optimal with respect to the entropy criterion.
\end{thm}

\section{Numerical experiments}
  \label{sec:sec5}

\begin{figure}[!t]
\begin{center}
\leavevmode
\hbox{
\epsfig{file=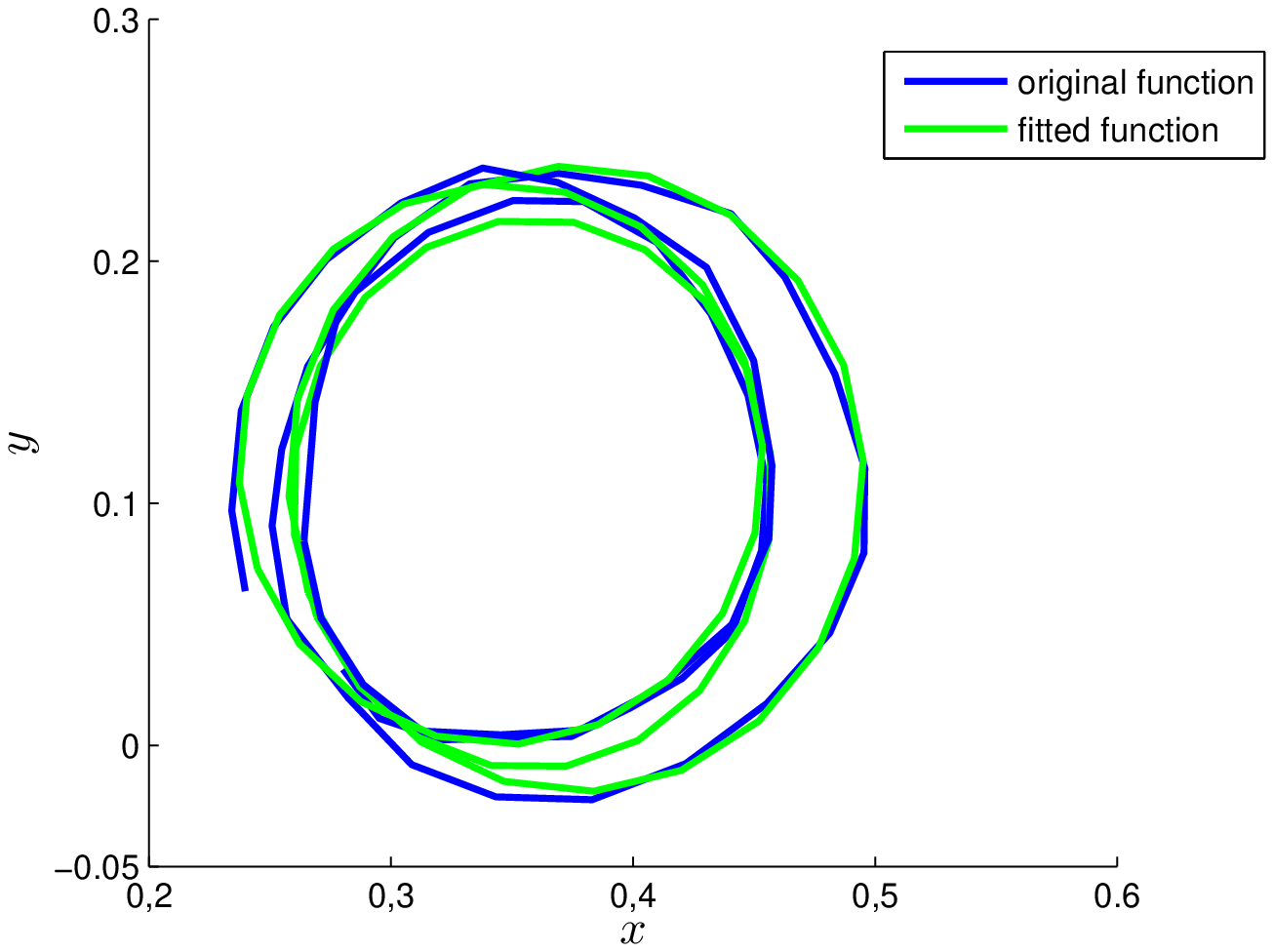,height=6cm, width=7.5 cm} \quad
\epsfig{file=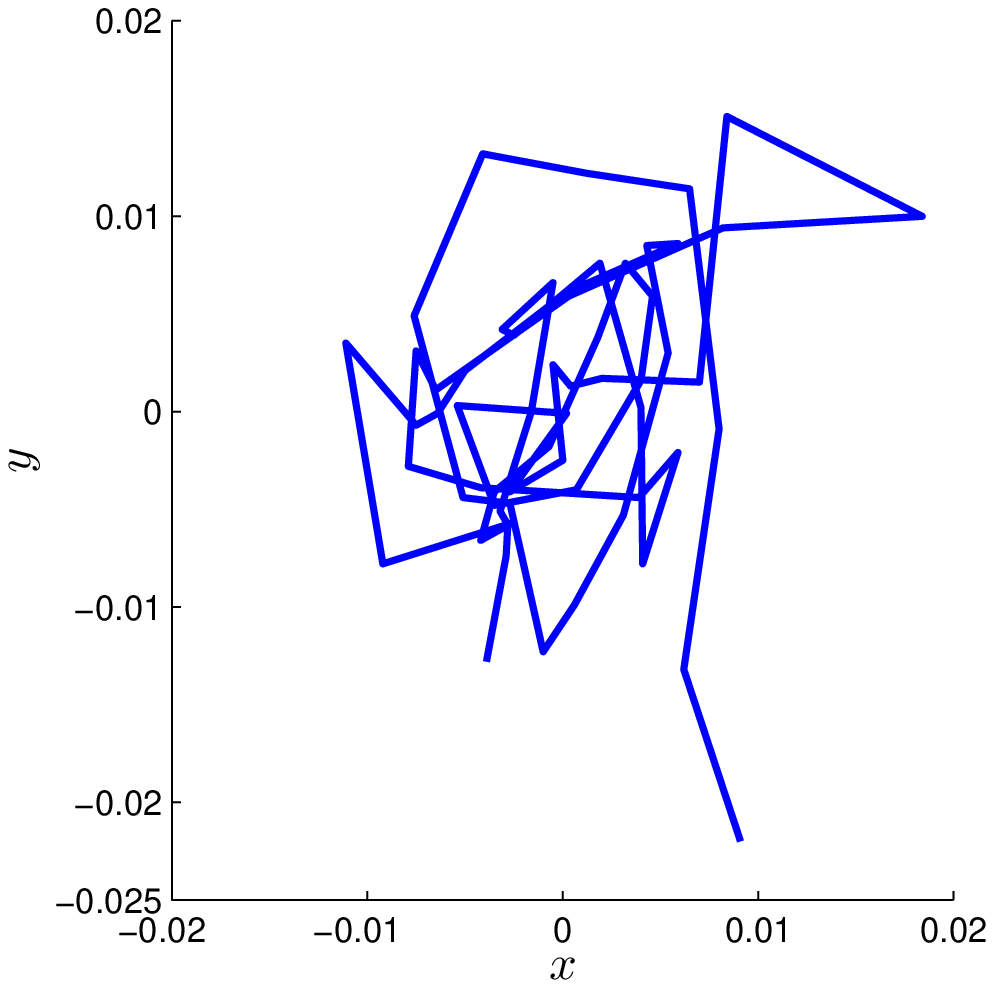,height=6cm}}

\centerline{\hbox to 9 truecm {\scriptsize (a) \hfill (b)}}

\end{center}
\caption{Real and fitted yearly polar motion (a) 2015-2017 based on data from IERS EOP C01 IAU2000, together with
the corresponding deviation plot (b).}
\label{fig2}
\end{figure}

  \begin{table}[!t]
\begin{center}
{\footnotesize
\begin{tabular}{|c|l|c|c|c|}
\hline
\multicolumn{2}{|c|}{}&$\lambda=2.4522, \ \omega=-4.1274$&$\lambda=4.9968, \ \omega=-0.3561$&
$\lambda=4.9366, \ \omega=-5.7767$ \\
\multicolumn{2}{|c|}{}&(2017)&(2016)&(2015)\\\hline
&optimal&11416&33633&18388\\
$n=3$&equispaced&11416&33633&18388\\
&rel. eff. (\%)&100&100&100\\\hline
&optimal&14724&25472&16959\\
$n=4$&equispaced&15470&25473&17977\\
&rel. eff. (\%)&95.18&99.99&94.34\\\hline
&optimal&14152&16305&11785\\
$n=5$&equispaced&20226&16320&11785\\
&rel. eff. (\%)&69.97&99.91&89.15\\\hline
\end{tabular}
}
\end{center}
\caption{IMSPE values corresponding to the optimal and to the
  equispaced design and relative efficiency of the equispaced
  design.}
\label{tab1}
\end{table}
In order to compare the performances of the two examined criteria we compare the optimal values of 
\ $\imspe ({\widehat Z})$ \ calculated using
{\tt fmincon} function of Matlab to its values
corresponding to the equispaced design which is optimal for the
entropy criterion. In Table \ref{tab1} the values of IMSPE are given
for both designs together with the relative efficiency of the
equispaced monotonic design with respect to the optimal one for various sample
sizes and combinations of damping parameter \ $\lambda$ \ and frequency \ $\omega$.
Observing that for the three-point design ($n=3$), the efficiency of the equispaced design is 100\%, so
in this case the equidistant design is optimal. In the other cases it is noticeable, that 
the optimal values of \ $\imspe ({\widehat Z})$ \ are better than the optimal values of entropy criterion.
The damping parameter \ $\lambda$ \ and frequency \ $\omega$ used for the simulations are
estimated based on public pole coordinates from the IERS EOP C01 IAU2000. For the estimation as a first 
step we give the model that fits best the real polar motion in the least squares sense. Then, according
to the Kolmogorov's model \eqref{chandler} of the Chandler wobble, deviation of the real and fitted yearly polar motion 
results in the required complex OU process. The estimation of parameters based on the given complex OU process and
the necessary maximum likelihood estimators for \ $\lambda$ \  and  \ $\omega$ are given in \citet{arato2,arato}.
As an example Figures \ref{fig2} shows the real and fitted yearly polar motion based on data from the past three years, together with
the corresponding deviation modelled as a complex OU process.

\section{Conclusions}
 \label{sec:sec6}
We derive the exact form of the IMSPE for a shifted complex OU process on
a compact interval and show that optimal design for prediction based on
IMSPE may well differ from the equidistant one. This is in contrast
both to the D-optimal design for estimation \citep{bszs}
and to the case of the classical real OU process \citep{baz}, but
similar to the optimal design for the prediction of OU sheets on a 
monotonic set \citep{bss13}. We also investigate the
properties of the optimal design with respect to entropy criterion and we 
show that the optimal designs are equidistant. Simulations illustrate
selected cases of optimal designs for small number of sampling
locations. The damping parameter \ $\lambda$ \ and frequency parameter 
\ $\omega$ \ used for the simulations are estimated based on real
data (pole coordinates from the IERS EOP C01 IAU2000), which is a well known application of the complex OU process, namely Kolmogorov's model \eqref{chandler} of the Chandler wobble. 
Since the above discussed designs depend on values of
damping and frequency parameters, obtained optimal designs are only locally
optimal. Such knowledge may be crucial for experiments to
increase efficiency of design in practical setups. 

\medskip
\noindent
{\bf Acknowledgments.} \  \  This research has been supported by the Hungarian -- Austrian intergovernmental S\&T cooperation program T\'ET\_{}15-1-2016-0046.
Kinga Sikolya acknowledges the support of the \'UNKP-17-4 New National Excellence Program of the Ministry of Human Capacities. The project has been supported by the European Union, co-financed by the European Social Fund.
S\'andor Baran is grateful for the support of the J\'anos Bolyai Research Scholarship of the Hungarian Academy of Sciences.

\begin{appendix} 
\section{Appendix}
  \label{sec:secA}
\subsection{Proof of Theorem \ref{IMSPE}}
To shorten our formulae, in what follows instead of
$Q(x,t_i)$ we are using simply $Q_i, \ i=1,2,\ldots ,
n$. According to the results of \citet{bszs}
\begin{equation*}
C(n)=
     \begin{bmatrix}
        {\mathbb I}_2 & {\mathrm e}^{A^{\top}\!d_1} & {\mathrm e}^{A^{\top}\!(d_1+d_2)} & 
        {\mathrm e}^{A^{\top}\!(d_1+d_2+d_3)} &\dots &\dots & 
        {\mathrm e}^{A^{\top}\!(\sum_{j=1}^{n-1}d_j)} \\
        {\mathrm e}^{Ad_1} &{\mathbb I}_2 &{\mathrm e}^{A^{\top}\!d_2} &
        {\mathrm e}^{A^{\top}\!(d_2+d_3)} &\dots &\dots &{\mathrm e}^{A^{\top}\!(\sum_{j=2}^{n-1}d_j)} \\
        {\mathrm e}^{A(d_1+d_2)} &{\mathrm e}^{Ad_2} &{\mathbb I}_2 
        &{\mathrm e}^{A^{\top}\!d_3} &\dots &\dots &{\mathrm e}^{A^{\top}\!(\sum_{j=3}^{n-1}d_j)} \\
         {\mathrm e}^{A(d_1+d_2+d_3)} & {\mathrm e}^{A(d_2+d_3)} 
         &{\mathrm e}^{Ad_3} &{\mathbb I}_2 &\dots &\dots &\vdots \\
        \vdots &\vdots &\vdots &\vdots &\ddots & &\vdots \\
        \vdots &\vdots &\vdots &\vdots & &\ddots &{\mathrm e}^{A^{\top}\!d_{n-1}} \\
        {\mathrm e}^{A(\sum_{j=1}^{n-1}d_j)} &
        {\mathrm e}^{A(\sum_{j=2}^{n-1}d_j)}
        &{\mathrm e}^{A(\sum_{j=3}^{n-1}d_j)} &\dots &\dots &
        {\mathrm e}^{Ad_{n-1}} &{\mathbb I}_2
     \end{bmatrix}.
\end{equation*}
and the inverse of \ $C(n)$ \ is given by
\begin{equation*}
  C^{-1}(n)=
     \begin{bmatrix}
        U_1 & -{\mathrm e}^{A^{\top}\!d_1}U_1 &0 & 0&\dots &\dots &0 \\
         -{\mathrm e}^{Ad_1}U_1 &V_2 &-{\mathrm e}^{A^{\top}\!d_2}U_2 &0 &\dots &\dots &0 \\
        0 & -{\mathrm e}^{Ad_2}U_2 &V_3  &-{\mathrm e}^{A^{\top}\!d_3}U_3&\dots &\dots &0 \\
        0 &0 &-{\mathrm e}^{Ad_3}U_3 &V_4&\dots &\dots &\vdots \\
        \vdots &\vdots &\vdots&\vdots  &\ddots & &\vdots \\
        \vdots &\vdots &\vdots &\vdots & &V_{n-1} &-{\mathrm e}^{A^{\top}\!d_{n-1}}U_{n-1} \\
        0 &0 &0 &\dots&\dots &-{\mathrm e}^{Ad_{n-1}}U_{n-1} &U_{n-1}
     \end{bmatrix},
\end{equation*}
where \ $U_k:=\big[{\mathbb I}_2-{\mathrm e}^{(A+A^{\top}\!)d_k}\big]^{-1},\ k=1,2,\dots,n-1$, \  and \ $V_k:=U_k+{\mathrm e}^{(A+A^{\top}\!)d_{k-1}}U_{k-1}, \ k=2,3,\dots,n-1$.

Consider first the $\mspe\big (\widehat Z(x)\big)$ given by
\eqref{mspe}. Short matrix algebraic calculations show
\begin{align*}
\left[
      \begin{BMAT}{c.c}{c.c}
      {\mathbb O}_2 & H_n\\
       H_n^{\top} & C(n)
      \end{BMAT}
\right]^{-1}\!\!\!\!\!&=\!
\left[
      \begin{BMAT}{c.c}{c.c}
      {\mathbb O}_2 & {\mathbb O}_{2n \times 2}^{\top}\\
      {\mathbb O}_{2n \times 2} & C^{-1}(n)
      \end{BMAT}
\right]-\big(H(n)C(n)^{-1}H(n)^{\top}\big)^{-1}\!\! \\
&\phantom{=========}\times \left[
      \begin{BMAT}{c.c}{c.c}
      {\mathbb I}_2 & -H(n)C^{-1}(n)\\
      -C^{-1}(n)H(n)^{\top} & C^{-1}(n)H(n)^{\top}H(n)C^{-1}(n)
      \end{BMAT}
\right],
\end{align*}
and according to the results of \citet{bszs}
we have
\begin{equation*}
H(n)C(n)^{-1}H(n)^{\top}\!\!\!=G(n){\mathbb I}_2\!=\!\Big(1\!+\!\sum_{\ell=1}^{n-1} g(d_{\ell})\Big){\mathbb I}_2, \quad \text{where}  \quad g(x)\!:=\!\frac{1\!-\!2{\mathrm e}^{-\lambda x}\cos(\omega x)\!+\!{\mathrm e}^{-2\lambda x}}{1\!-\!{\mathrm e}^{-2\lambda x}}, 
\end{equation*}
$x>0, \ \ \text{and} \ \ g(0)\!:=\!0$.\ We remark that
$H(n)C(n)^{-1}H(n)^{\top}$ is the FIM on parameter vector \ $(m_1,m_2)^{\top}$ \ based on observations \ $\big\{ \big(Z_1(t_j),Z_2(t_j)\big), \ j=1,2,\ldots ,n\big\}$. \ In this way, we obtain
\begin{align*}
\mspe \big (\widehat Z(x)\big)=&\,\frac{\sigma^2}{2\lambda}{\mathrm{tr}}\Bigg[{\mathbb I}_2\!-\!{\mathcal
  Q}(x)C^{-1}(n) {\mathcal Q}^{\top}(x)\!+\! G(n)^{-1}\bigg({\mathbb I}_2-H(n)C^{-1}(n){\mathcal Q}^{\top}(x) \\
& -{\mathcal Q}(x)C^{-1}(n)H(n)^{\top}
+{\mathcal Q}(x)C^{-1}(n)H(n)^{\top}H(n)C^{-1}(n){\mathcal Q}^{\top}(x)\bigg)\Bigg] \\
=&\,\frac{\sigma^2}{2\lambda}{\mathrm{tr}}\Bigg[{\mathbb I}_2\!-\!{\mathcal
  Q}(x)C^{-1}(n) {\mathcal Q}^{\top}(x)\!+\! G(n)^{-1}\bigg({\mathbb I}_2-{\mathcal Q}(x)C^{-1}(n)H(n)^{\top}\bigg) \\
& \times \bigg({\mathbb I}_2-{\mathcal Q}(x)C^{-1}(n)H(n)^{\top}\bigg)^{\top}\Bigg].
\end{align*}	
Examining member by member the above expression
\begin{align*}
\tr\Big(&{\mathcal Q}(x)C^{-1}(n){\mathcal Q}^{\top}(x)\Big)= \\
=&\,\tr\Bigg({Q_n}{Q^{\top}_n}+\sum_{i=1}^{n-1}U_i\Big({Q_i}{Q^{\top}_i}-
{Q_i}{\mathrm e}^{A^{\top}\!d_i}{Q^{\top}_{i+1}}-{Q_{i+1}}{\mathrm e}^{Ad_i}{Q^{\top}_{i}}
+{Q_{i+1}}{\mathrm e}^{(A+A^{\top})d_i}{Q^{\top}_{i+1}}\Big)\Bigg)\\
=&\,2p_{n,n}+2\sum_{i=1}^{n-1}\frac
{p_{i,i}-2{\pi_i}p_{i+1,i}+{{\pi_i}^2}p_{i+1,i+1}}{1-{\pi_i}^2}, \\[3mm]
  \end{align*}
\begin{align*}
\tr\Big(&H(n)C^{-1}(n){\mathcal Q}^{\top}(x)+{\mathcal Q}(x)C^{-1}(n)H(n)^{\top}\Big)=\\
=&\,\tr\Bigg(\!{Q_n}\!+\!{Q^{\top}_n}\!+\!\sum_{i=1}^{n-1}(Q_i-Q_{i+1}{\mathrm e}^{Ad_i})U_i({\mathbb I}_2-{\mathrm e}^{A^{\top}\!d_i})
\!+\!\sum_{i=1}^{n-1}U_i({\mathbb I}_2-{\mathrm e}^{Ad_i})(Q^{\top}_i-{\mathrm e}^{A^{\top}\!d_i}Q^{\top}_{i+1})\!\Bigg)\\
  =&\,4q_{n,n}+4\sum_{i=1}^{n-1}\frac{(q_{i,i}-{\pi_i}q_{i,i+1})-{\pi_i}(q_{i+1,i}-{\pi_i}q_{i+1,i+1})}{1-{\pi_i}^2}, \\[3mm]
\tr\Big(&{\mathcal Q}(x)C^{-1}(n)H(n)^{\top}H(n)C^{-1}(n){\mathcal Q}^{\top}(x)\Big)=\\
=&\,\tr\Bigg(\bigg(\sum_{i=1}^{n-1}(Q_i\!-\!Q_{i+1}{\mathrm e}^{Ad_i})U_i({\mathbb I}_2\!-\!{\mathrm e}^{A^{\top}\!d_i})\!+\!Q_{n}\bigg)
   \bigg({Q^{\top}_n}\!+\!\sum_{i=1}^{n-1}U_i({\mathbb I}_2\!-\!{\mathrm e}^{Ad_i})(Q^{\top}_i\!-\!{\mathrm e}^{A^{\top}\!d_i}Q^{\top}_{i+1})\bigg)\Bigg)\\
=&\,2p_{n,n}+4\sum_{i=1}^{n-1}\frac{\big(p_{n,i}-\pi_ip_{n,i+1}\big)\big(\cos(\omega(t_n-t_i))-{\pi_i}\cos(\omega(t_n-t_{i+1}))\big)}{1-{\pi_i}^2}\\
&+2\sum_{i=1}^{n-1} \frac{\big(p_{i,i}-2{\pi_i}p_{i+1,i} +{\pi_i}^2p_{i+1,i+1}\big)\big(1-2{\pi_i}\cos(\omega d_i)+{\pi_i}^2\big)}{(1-\pi_i^2)^2} \\
&+4\sum_{i=2}^{n-1}\sum_{j=1}^{i-1} \big(p_{i,j}-\pi_ip_{i+1,j}-\pi_jp_{i,j+1}+\pi_i\pi_j p_{i+1,j+1}\big)\\
&\phantom{+4\sum_{i=2}^{n-1}}\times \frac{\big(\cos(\omega(t_{i}\!-\!t_j))\!-\!{\pi_i}\cos(\omega(t_{i+1}\!-\!t_j))\!-\!{\pi_j}\cos(\omega(t_i\!-\!t_{j+1}))\!+\!{\pi_i}{\pi_j}\cos(\omega (t_{i+1}\!-\!t_{j+1}))}{(1-{\pi_i}^2)(1-{\pi_j}^2)}.
\end{align*}
where \ $\pi_i:=\exp(-\lambda d_i)$, \ with \ $d_i:=t_{i+1}-t_i$, \ $p_{i,j}:=\exp(-\lambda(|x-t_i|+|x-t_j|))$, \ and \ $q_{i,j}:=\exp(-\lambda|x-t_i|)\cos(\omega(x-t_j))$, \
 \ $i,j=1,2,\dots,n-1$.
Further, according to the definition of the IMSPE criterion, we can write
\begin{equation*}
 \imspe\big (\widehat Z\big)=2\big(1-A_n+G(n)^{-1}B_n\big),
\end{equation*}
where
\begin{align*}
A_n:=&\, \frac{1}{2}\int\limits_{\mathcal X}\tr\Big({\mathcal
  Q}^{\top}(x)C^{-1}(n){\mathcal Q}(x)\Big){\mathrm d}x, \\
B_n:=&\, \frac{1}{2}\int\limits_{\mathcal X}\tr\bigg({\mathbb I}_2-H(n)C^{-1}(n){\mathcal Q}(x) 
 -{\mathcal Q}^{\top}(x)C^{-1}(n)H(n)^{\top} \\
&\phantom{=======}+{\mathcal Q}^{\top}(x)C^{-1}(n)H(n)^{\top}H(n)C^{-1}(n){\mathcal Q}(x)\bigg){\mathrm d}x.
\end{align*}
In this way, using that for $i\land j:=\min\{i,j\}, \ i\lor j:=\max\{i,j\}, \ i,j
\in{\mathbb N}$, \ we have
\begin{align*}
\int\limits_{\mathcal
  X}\!p_{i,j}\, {\mathrm d}x \!= &\int\limits_0^1\!
  {\mathrm e}^{-\lambda (|x-t_i|+|x-t_j|)} {\mathrm d}x \\
=&\, \frac{1}{2\lambda} \Big(2{\mathrm e}^{-\lambda (t_i-t_j)} - {\mathrm e}^{-\lambda(t_i+t_j)}- {\mathrm e}^{-\lambda(1-t_i+1-t_j)}\Big) 
+\big( t_i- t_j\big){\mathrm e}^{-\lambda(t_i-t_j)} \\ =&\, \frac{1}{2\lambda} \Big(2{\mathrm e}^{-\lambda (d_j+ \cdots + d_{i-1})} - {\mathrm e}^{-\lambda(2d_1+\cdots +2d_{j-1}+d_j+\cdots +d_{i-1})}- {\mathrm e}^{-\lambda(d_j+\cdots +d_{i-1}+2d_i+\cdots +2d_{n-1})}\Big) \\
                                 &+\big(d_j+ \cdots + d_{i-1}\big){\mathrm e}^{-\lambda(d_j+ \cdots + d_{i-1})}:=\varrho_{i,j},
\quad  \quad 1\leq j \leq i \leq n,  
\end{align*}
\begin{align*}
\int\limits_{\mathcal
  X}\!q_{i,j}\, {\mathrm d}x \!=& \int\limits_0^1\!
  {\mathrm e}^{-\lambda |x-t_i|}\cos(\omega (x-t_j)) {\mathrm d}x=\frac {2\lambda}{\lambda^2 +\omega^2}\cos(w(t_i-t_j)) \\
&+\!\frac {{\mathrm e}^{-\lambda t_i}}{\lambda^2 +\omega^2}\big(\omega \sin(\omega t_j)\!-\!\lambda \cos (\omega t_j)\big)\!+\!\frac{{\mathrm e}^{-\lambda (1-t_i)}}{\lambda^2 +\omega^2}\big(\omega \sin(\omega (1-t_j))\!-\!\lambda \cos (\omega (1-t_j))\big)\\
=&\,\frac {2\lambda}{\lambda^2 +\omega^2}\cos(w(d_{i\land j}+\ldots +d_{i \lor j-1})) \\
&+\frac {{\mathrm e}^{-\lambda (d_1+ \cdots + d_{i-1})}}{\lambda^2 +\omega^2}\big(\omega \sin(\omega (d_1+ \cdots + d_{j-1}))-\lambda \cos (\omega (d_1+ \cdots + d_{j-1}))\big) \\
&+\frac {{\mathrm e}^{-\lambda (d_i+ \cdots + d_{n-1})}}{\lambda^2 +\omega^2}\big(\omega \sin(\omega (d_j+ \cdots + d_{n-1}))-\lambda \cos (\omega (d_j+ \cdots + d_{n-1}))\big):=v_{i,j},
\end{align*}	
with the empty sum to be defined as zero, short calculation leads to \eqref{an} and \eqref{bn}.
\proofend

\subsection{Proof of Theorem \ref{entropy}}

Following the idea of \citet[Lemma 3.1]{baz} the covariance matrix \ $C(n)$ \ can be written as
\begin{equation*}
C(n)=LDU,
\end{equation*}
where \ $L$ \ is a lower, \ $U$ \ is an upper block triangular matrix with \ ${\mathbb I}_2$ \ as blocks in the main diagonal and \ $D$ \ is block diagonal matrix with
\begin{equation*}
\bigg({\mathbb I}_2, {\mathbb I}_2-{\mathrm e}^{(A+A^{\top}\!)d_1}, 
{\mathbb I}_2-{\mathrm e}^{(A+A^{\top}\!)d_2}, \ldots , 
{\mathbb I}_2-{\mathrm e}^{(A+A^{\top}\!)d_{n-1}}\bigg)
\end{equation*}
as main block matrices.
In this way 
\begin{equation*}
\det C(n) =\det {\mathbb I}_2 \times \det \big({\mathbb I}_2-{\mathrm e}^{(A+A^{\top}\!)d_1}\big) \times \ldots 
\det \big({\mathbb I}_2-{\mathrm e}^{(A+A^{\top}\!)d_{n-1}}\big)  =\prod_{i=1}^{n-1} \big(1-2\pi_i^2\big),
\end{equation*}
which proves \eqref{entropy}.
Now, from \eqref{entropy} we have
\begin{equation*}
\ln \det C(n)=\sum_{i=1}^{n-1} f(d_i), \qquad \text{where} \quad f(x):=\ln \big(1-2{\mathrm e}^{-2\lambda x}\big).
\end{equation*}
Since
\begin{equation*}                                      
\frac {\partial^2 f(d)}{\partial d^2}=-\frac {8\lambda^2 {\mathrm e}^{2\lambda d}}{({\mathrm e}^{2\lambda d}-2)^2}<0 
 \quad \text{for any} \quad d\in(0,1),
\end{equation*}
$f(x)$ \ is a concave function of \ $x$, \ and the result follows from Schur-concavity of the entropy criterion.
\proofend

\end{appendix}

\end{document}